\documentclass{article}
\usepackage{amsthm}
\usepackage{amsfonts}
\usepackage{amssymb}
\usepackage{amsgen}
\usepackage{amsmath}
\usepackage{amsopn}
\usepackage{verbatim}
\usepackage{xypic}
\usepackage{xspace}
\usepackage{multicol}
\usepackage{makeidx}
\usepackage{eepic}
\usepackage{upref}
\usepackage{url}
\usepackage[mathscr]{eucal}

\theoremstyle{plain}
\newtheorem{thm}{Theorem}

\newtheorem{prop}[thm]{Proposition}

\newtheorem{conj}[thm]{Conjecture}
\theoremstyle{definition}

\newtheorem{defn}[thm]{Definition}

\newtheorem{eg}[thm]{Example}

\newcommand{\tn}{{\tilde{n}}}

\newcommand{\ors}{{\vec{s}}}

\newcommand{\N}{{\mathbbN}}

\DeclareMathOperator{\wt}{{wt}}

\newcommand{\gt}{{\tau}}

\renewcommand{\N}{{\mathbb N}}

\newcommand{\Lra}{\Longrightarrow}

\begin{document}


\begin{center}
 {\huge Examples of finite $p$-divisible sets of MHS}
\end{center}

\begin{center}
 {\large Department of Mathematics, Eckerd College, St. Petersburg, FL 33711, USA}
\end{center}

\vskip0.6cm
\begin{center}
 {\large Jianqiang Zhao}
\end{center}

\vskip1cm The following examples give evidence to the following
conjecture contained in my paper \cite{main}. All the main
theoretical results can be found in that paper.

\begin{conj}
Let $d$ be a positive integer and $\ors \in {\mathbb N}^d$. Then
the set $J(\ors|p)$ is finite for every prime $p$.
\end{conj}

\begin{eg} \label{eg:zeta2}
The first example we would like to do is about the partial
sums of $\zeta(2)$ series. The prime $p=7$ is a little
different from the others because $H(2;3)=49/36$. We find that
$$J_1(2|7)=\{3,6\}, J_1(2|7^2)=\{3\}.$$
By \cite[Thm~3.6]{main} $n=7\tn+r\in J_2(2|7)$, $0\le r\le 6$,
exists if and only if $\tn=3,$ and
$$H(r)\equiv -\psi(3)\equiv -1/36 \equiv -1 \pmod{7}.$$
 From the congruence
$$H(2;1)\equiv 1, \quad H(2;2)\equiv 3,\quad
H(2;4)\equiv 4, \quad H(2;5)\equiv -1 \pmod{7},$$ we find $r=5$.
But it's an easy matter to check that
$$H(2;26) \equiv 14 \not \equiv 0\pmod{7^2}.$$
Hence $J_2(2|7^2)=\emptyset$ and consequently
$J(2|7)=\{0,3,6,26\}$.

For all the other primes $p\ne 7$ from 5 up to 1061\footnote{We choose 1061
because it's the first irregular prime greater than 1000.} we find that
$$J_1(2|p)=\{0,(p-1)/2, p-1\}\cup T(2|p)\cup \{p-1-r: r\in T(2|p)\}$$
where $T(2|p)$
are listed in Table~\ref{Ta:J2p} if $T(2|p)\ne \emptyset$.
Moreover, $J_1(2|p^2)=\emptyset$ which implies
$J(2|p)=\{(p-1)/2, p-1\}\cup T(2|p)\cup \{p-1-r: r\in T(2|p)\}$ in
this range.

\begin{table}[h]
\begin{center}
\begin{tabular}{  ||c|c||c|c||c|c||c|c||c|c||c|c|| }
 \hline
$p$&$T(2|p)$&$p$&$T(2|p)$&$p$& $T(2|p)$&$p$&$T(2|p)$&$p$&$T(2|p)$&$p$& $T(2|p)$\\ \hline
37&15      &163&61   &419&111 &563&175,227 &677&153 &883&151  \\ \hline
41&4       &167&61    &421&59   &569&199    &709&123  &911&345  \\ \hline
43&11      &181&85    &433&179  &571&247    &727&197,239&929&64 \\ \hline
59&6,24    &211&99    &457&216 &577&134,158 &739&93     &953&199 \\ \hline
97&15      &241&60,96 &467&158,170&601&17   &787&344    &967&463 \\ \hline
107&39     &269&50    &479&5      &617&97   &797&185,226&971&429 \\ \hline
127&23     &307&27    &487&32   &619&16,70,286 &811&371 &991&72  \\ \hline
137&44     &311&43    &491&173  &643&17,222    &821&39  &997&205,310 \\ \hline
149&37     &373&54    &499&134  &653&246,307   &859&414 &1013&430 \\ \hline
157&25     &383&150   &547&165  &659&232       &863&226 &1031&384 \\ \hline
\end{tabular}
\caption{The set $J(2|p)$, $p\le 1061$, is finite.}
\label{Ta:J2p}
\end{center}
\end{table}
\end{eg}

\begin{eg} \label{eg:zeta3}
In this example we will determine $J(3|p)$ for primes
$p$ no greater than 1061. First let $p=3$.
Easy computation by hand 
will show that $J_1(3|3)=\{2\}$
and $J_1(3|3^3)=\emptyset$. Thus $J(3|3)=\{0,2\}$.
Let $p=5.$ We see that $J_1(3|5)=\{4\}$ and $J_1(3|5^3)=\emptyset$.
Hence $J(3|5)=\{0,4\}$.

Similarly, for $p$ between $7$ and $1061$ we get that
$J(3|p)=\{0,p-1\}\cup T(3|p)\cup \{p-1-r: r\in T(3|p)\}$ where
$T(3|p)$ are listed in Table~\ref{Ta:J3p} if $T(3|p)\ne
\emptyset$. In fact we find that when $p$ lies in this range
$J(3|p^2)=\{0,p-1\}$ exactly, and $J(3|p^3)=\{0\}$ except for
$J(3|37^3)=\{0,36\}$.  But by our \cite[Thm~2.6]{1stpart} we know
that there might be infinitely many other primes with
$J(3|p^3)=\{0,p-1\}$, as long as $(p, p-5)$ is an irregular pair.

Let's look at $p=37.$ Maple tells us that $J_1(3|37)=\{4, 13,23,
32, 36\}$, $J_1(3|37^2)=J_1(3|37^3)=\{36\}.$ By
\cite[Thm~3.6]{main}, in order for $n=37 \tn+r\in J_2(3|37)$ to
exist it is necessary and sufficient that
$$\tn=36 \mbox{ and } H(3,r)\equiv -\psi(\tn) \equiv 28 \pmod{37}.$$
A quick Maple computation picks up $r=8$ and $r=28$. (Recall
that all such $r$ should be paired and add up to $p-1=36$).
So $J_2(3|37)=\{37\cdot 36+8, 37\cdot 36+28\}=\{1340, 1360\}.$
Direct computation shows that these two number do not
lie in $J_2(3|37^3)$ which shows that $J_3(3|37)=\emptyset.$
To summarize we get
\begin{equation}\label{J3|37}
J(3|37)=\{0,4, 13, 23, 32, 36, 1340, 1360\}.
\end{equation}

\begin{table}[h]
\begin{center}
\begin{tabular}{  ||c|c||c|c||c|c||c|c||c|c||c|c|| }
 \hline
$p$&$T(3|p)$&$p$&$T(3|p)$&$p$&$T(3|p)$&$p$&$T(3|p)$&$p$&$T(3|p)$&$p$&$T(3|p)$\\ \hline
11&4    &163&73,76 &277&56    &421&76    &619&138  &863&135\\ \hline
17&7    &173&56    &281&45    &431&11,25 &673&147  &877&73    \\ \hline
31&8    &179&24    &283&113   &443&105   &677&324  &883&195   \\ \hline
37&see \eqref{J3|37}&181&76&317&103  &479&127   &743&307  &887&145   \\ \hline
47&5    &193&38    &331&154,161&521&43   &751&52   &911&287 \\ \hline
53&6    &223&7   &347&16,73,107&523&93   &773&311  &941&241,419\\ \hline
67&28,30&227&91    &349&83    &563&62,156&787&9,199&953&29,69 \\ \hline
89&10,43&251&3     &359&10,136&571&274   &809&213  &1013&345\\ \hline
113&39,43&263&76   &367&163   &577&76    &823&37    &1021&383   \\ \hline
137&44  &269&131   &379&160   &599&224  &829&256  &1051&284   \\ \hline
149&35  &271&70,105&389&63    &607&42   &839&71   &1061&515 \\ \hline
\end{tabular}
\caption{The set $J(3|p)=\{0,p-1\}\cup T(3|p)\cup
\{p-1-r: r\in T(3|p)\}$, $p\le 1061$, is finite.}
\label{Ta:J3p}

\end{center}
\end{table}

\end{eg}


\begin{eg} \label{eg:zeta4}
The partial sums of $\zeta(4)$ series are very regular. When $p=3$ and
$p=5$, we see that
$$J_1(4|p)=\emptyset\Lra J(4|p)=\{0\}.$$

For all the primes from 7 up to less than 1100 we have the following:
$$J_1(4|p)=\{(p-1)/2, p-1\}\cup T(4|p)\cup \{p-1-r: r\in T(4|p)\}$$
where $T(4|p)$
are listed in Table~\ref{Ta:J4p} if $T(2|p)\ne \emptyset$.
Moreover, $J_1(4|p^4)=J_1(4|p^2)=\emptyset$ for all primes
in this range. As a consequence
$$J(4|p)=\{0,(p-1)/2, p-1\}\cup T(4|p)\cup \{p-1-r: r\in T(4|p)\}\quad
\forall p<1100.$$
\begin{table}[h]
\begin{center}
\begin{tabular}{  ||c|c||c|c||c|c||c|c||c|c||c|c|| }
 \hline
$p$&$T(4|p)$&$p$&$T(4|p)$&$p$&$T(4|p)$&$p$&$T(4|p)$&$p$&$T(4|p)$&$p$&$T(4|p)$\\ \hline
17&2    &157&61  &307&32   &487&125    &757&115  &941&341,355  \\ \hline
41&18   &163&10  &311&89 &499&97,137,232&773&354 &947&241  \\ \hline
59&15   &181&74  &317&11  &503&60      &787&356  &967&90139 \\ \hline
67&24   &191&33  &331&134 &523&186    &797&43    &997&196   \\ \hline
71&28   &199&3,63&337&43  &541&111    &809&374   &1009&383  \\ \hline
79&6    &227&43  &359&53  &617&51     &827&320   &1021&359   \\ \hline
97&38   &229&67  &401&188 &619&255    &829&84    &1031&141 \\ \hline
101&46  &239&42  &409&99  &647&40,194 &839&88,238&1033&22    \\ \hline
103&47  &251&16,88&431&13 &691&292    &863&252   &1039&107,391  \\ \hline
131&58  &257&112  &461&209&701&160    &883&212   &1049&404   \\ \hline
137&10,51&277&91  &463&30 &719&354    &919&250   &1061&212  \\ \hline
139&63  &283&11  &479&54  &743&281    &937&89    &1093&203 \\ \hline
\end{tabular}
\caption{The set $J(4|p)=\{0,(p-1)/2,p-1\}\cup T(4|p)\cup
    \{p-1-r: r\in T(4|p)\}$, $p\le 1100$, is finite.}
\label{Ta:J4p}
\end{center}
\end{table}
\end{eg}


\begin{eg}  \label{eg:zeta5}
Let's turn to the partial sums of $\zeta(5)$ series. First we easily
have $J_1(5|3)=\{2\}$ and $J_1(5|3^3) =\emptyset$ so that
$J(5|3)=\{0,2\}$. Furthermore
$$H(5;1)=1,\ H(5;2)\equiv 44,\ H(5;3)\equiv 26,\
H(5;4)\equiv 1 \pmod{5^4}.$$
So we get $J_1(5|5)=J_1(5|5^2)=J_1(5|5^3)=\{4\}$, but
$J_1(5|5^5)=\emptyset$. So $J(5|5)=\{0,4\}.$
Now we jump to $p=37$ because it gives the first irregular pair
$(37,37-5)$ and therefore shows a different pattern.
We must have $18\in J_1(5|37)$ as guaranteed by
Proposition~4.15(a). Computation shows that
in fact $J_1(5|37)=\{6,9,12,18,24,27,30,36\}.$
We also have $J(5|37^2)=\{0,6,36\}$ which is a little strange,
as the set is  obviously asymmetric. Finally $J_1(5|37^5)=\emptyset$
which implies that
$$J(5|37)=\{0,6,9,12,18,24,27,30,36\}.
$$

In general, when primes $p$ ranges from 7 to 1061 but
$p\ne 37$ we find by Maple
$$J(5|p)=\{0,p-1\}\cup T(5|p)\cup \{p-1-r: r\in T(5|p)\},
    \quad J(5|p^2)=\{0,p-1\},\quad J(5|p^3)=\emptyset,$$
where $T(5|p)$ is given by Table~\ref{Ta:J5p}
whenever $T(5|p)\ne \emptyset$.
The inclusion $p-1\in J(5|p^2)$ is enforced by
Lemma~2.2 and Theorem~2.6.

\begin{table}[h]
\begin{center}
\begin{tabular}{  ||c|c||c|c||c|c||c|c||c|c||c|c|| }
 \hline
$p$&$T(5|p)$&$p$&$T(5|p)$&$p$&$T(5|p)$&$p$&$T(5|p)$&$p$&$T(5|p)$&$p$&$T(5|p)$\\ \hline
11&2        &179&63  &317&142,154&467&179  &743&299,334&929&268   \\ \hline
29&6        &197&27   &349&111    &479&94,100&761&265   &937&349 \\ \hline
47&14       &223&32,35&359& 55    &563&42    &773&217   &947&26\\ \hline
71&9        &241&49   &379&187    &571&43    &811&61    &967&6,24  \\  \hline
83&3,15,21  &257&26   &397&48     &613&99    &821&180   &991&94    \\  \hline
97&3,22     &269&128  &409&151    &641&97    &853&168   &1009&66,451\\  \hline
107&37      &283&38   &419&95     &659&137,286&857&274 &1021&449 \\  \hline
127&19      &307&124  &431&151    &683&115   &859&89   &1031&139 \\  \hline
149&14      &311&144  &433&197     &691&273   &877&324  &1033&362 \\  \hline
157&32      &313&10,95 &439&211  &709&344   &883&54  &1039&181\\ \hline
\end{tabular}
\caption{The set $J(5|p)=\{0,p-1\}\cup T(5|p)\cup \{p-1-r: r\in
T(5|p)\}$, $p\le 1061$, is finite.}
\label{Ta:J5p}
\end{center}
\end{table}
\end{eg}

Similar but larger range of data is available online \cite{online}.
To summarize we have

\begin{prop} \label{prop:eg}
Let $p$ be a prime such that $p\le 3001$. Then $J(s|p)$ is
finite for $2\le s\le 300$.
\end{prop}

\pagebreak

For other types of $\ors$ we need the Criterion Theorem in our
main article \cite{main}.

\noindent
{\bf Criterion Theorem}. {\em Let $d\ge 2$ be a positive
integer and $p$ be a prime such that $d\in G_{t_0}$. Let
$\ors=(s_1,\cdots,s_d)\in \N^d$ and put $m=\min\{s_i: 1\le i\le
d\}$. For $t\in \N$ set
$f(\ors,p;t)=\min\{-v_p\big(H(\ors;n)\big): n\in G_t\}$. If there
is $\gt> t_0$ such that}
$$f(\ors,p;\gt)> (\wt(\ors)-m)(\gt-1) -m,$$
{\em then $J(\ors|p)$ is finite.}

\begin{prop} \label{prop:Jsdbig}
Let $s$ and $d$ be two positive integers. Suppose $2\le s\le 6$
and $d\le 5$. Then the $p$-divisible set $J(\{s\}^d|p)$ is finite
for the consecutive five primes immediately after $sd+2$. Moreover
there's always some prime $p$ such that
$J(\{s\}^d|p)=RJ(\{s\}^d;p)$ where
$$RJ(\{s\}^d;p)=\begin{cases}
\{0,p-1\} \quad &\text{if } 2\nmid s,\\
\{0,i+(p-1)/2,p-1,p: 0\le i\le d-1\} \quad &\text{if }2|s.
\end{cases}
$$
\end{prop}
\begin{proof}
We may assume $d\ge 2$ by Prop.~\ref{prop:eg}. Then the
proposition follows from Criterion Theorem by computer computation
whose results are listed in Tables~\ref{Ta:2dp}-{Ta:6dp}. We put
the Maple code at the end of this supplement. In each row of the
tables we list the number $\gt$ used when applying Criterion
Theorem and the discrepancy $J(\{s\}^d|p)\setminus RJ(\{s\}^d;p)$,
separated by a semicolon. We include the data for $s=1$ for future
reference where we list $RJ$ itself instead of the discrepancy when
$p\le d+2$.

\begin{table}[h]
\begin{center}
\begin{tabular}{  ||c|lc|lc|lc|lc||  } \hline
 $p $& \ &$d=2$  &\ &$d=3$   &\ &$d=4$   &\ &$d=5$  \\ \hline
 3& 6,&\{0,5\} &10,& \{0,8\} &$>$ 12, &\{0\}&$>$10, &\{0,8\} \\ \hline
 7& 4,&\{4\} &$>$8, &\{11\} & $>$8, & $\emptyset$ &$>8$, & \{0,6,14\} \\ \hline
 11&$>$8 &$\emptyset$&$>$8 &$\emptyset$ &$>$8, & $\emptyset$  &$>$8, & $\emptyset$  \\ \hline
 13& 4,&$\emptyset$&$>$8, &$\emptyset$& $>$8, & $\emptyset$ & $>$8, & $\emptyset$  \\ \hline
 17&$>$8 &\{11,13\}&$>$8, &\{5\}&$>$8,&$\emptyset$& $>$8, & $\emptyset$  \\ \hline
 31& 4,&\{17,22\} & \  & \  & \ & \  & \  & \ \\ \hline
\end{tabular} \caption{Finiteness of $J(\{1\}^d|p)$ for
$2\le d\le 5$.} \label{Ta:1dp}
\end{center}
\end{table}

\nopagebreak
\begin{table}[h]
\begin{center}
\begin{tabular}{  ||c|lc|lc|lc|lc||  } \hline
 $p $& \ &$d=2$  &\ &$d=3$   &\ &$d=4$   &\ &$d=5$  \\ \hline
 7& 3,&$\emptyset$ & \  & \  & \ & \  & \  & \ \\ \hline
 11&2,&$\emptyset$&3,&$\emptyset$&4,&$\emptyset$&\ &\ \\ \hline
 13&2,&$\{4\}$&3,&$\emptyset$&4,&$\emptyset$&5,&$\emptyset$\\ \hline
 17&2,&$\emptyset$&3,&$\{14\}$&4,&$\emptyset$&5,&$\emptyset$\\ \hline
 19&2,&$\emptyset$&3,&$\{14\}$&4,&$\emptyset$&5,&$\emptyset$\\ \hline
 23&2,&$\{17\}$&3,&$\emptyset$&4,&$\emptyset$&5,&$\emptyset$\\ \hline
 29&2,&$\emptyset$&3,&$\emptyset$&4,&$\{11\}$&5,&$\emptyset$\\ \hline
 31&2,&$\emptyset$&4,&\{5,11\}&\ &\ &\ &\ \\ \hline
 37&2,&$\{6\}$&\ &\ &\ &\ &\ &\ \\ \hline
 41&2,&$\emptyset$&\ &\ &\ &\ &\ &\ \\ \hline
\end{tabular} \caption{Finiteness of $J(\{2\}^d|p)$ for
$2\le d\le 5$.} \label{Ta:2dp}
\end{center}
\end{table}

\begin{table}[h]
\begin{center}
\begin{tabular}{  ||c|lc|lc|lc|lc||  } \hline
 $p $& \ &$d=2$  &\ &$d=3$   &\ &$d=4$   &\ &$d=5$  \\ \hline
 11&2,&$\emptyset$&\ &\ &\ &\ &\ &\ \\ \hline
 13&2,&$\emptyset$&3,&\{10\}&\ &\ &\ &\ \\ \hline
 17&2,&\{6\}&3,&\{10,13\}&4,&$\emptyset$&\ &\ \\ \hline
 19&2,&\{12\}&3,&$\emptyset$&4,&$\emptyset$&5,&\{10,16\}\\ \hline
 23&2,&$\emptyset$&3,&\{7,15\}&4,&$\emptyset$&5,&\{14\}\\ \hline
 29&2,&$\emptyset$&3,&$\emptyset$&4,&\{6,9,23\}&5,&\{11\}\\ \hline
 31&2,&\{11\}&\ &\ &4,&\{18,28\}&5,&$\emptyset$\\ \hline
 37&3,&\{13,14,21,23,24,28,73\} &\ &\ &\ &\ &5,&$\emptyset$\\ \hline
\end{tabular} \caption{Finiteness of $J(\{3\}^d|p)$ for
$2\le d\le 5$.} \label{Ta:3dp}
\end{center}
\end{table}

\begin{table}[h]
\begin{center}
\begin{tabular}{  ||c|lc|lc|lc|lc||  } \hline
 $p $& \ &$d=2$  &\ &$d=3$   &\ &$d=4$   &\ &$d=5$  \\ \hline
 11&2,&$\emptyset$&\ &\ &\ &\ &\ &\ \\ \hline
 13&2,&$\{4,10\}$&\ &\ &\ &\ &\ &\ \\ \hline
 17&2,&\{13\}&3,&\{6,13,14\}&\ &\ &\ &\ \\ \hline
 19&2,&$\emptyset$&3,&$\emptyset$&4,&$\emptyset$&\ &\ \\ \hline
 23&2,&$\emptyset$&3,&$\emptyset$&4,&$\emptyset$&5,&$\emptyset$\\ \hline
 29&2,&$\emptyset$&3,&$\emptyset$&4,&\{24\}&5,&\{11,24,25\}\\ \hline
 31&2,&$\emptyset$&3,&\{21,24\}  &4,&\{28\}&5,&$\emptyset$\\ \hline
 37&2,&$\emptyset$&\ &\ &4,&$\emptyset$ &5,&\{31\}\\ \hline
 41&2,&\{40\}&\ &\ &\ &\ &5,&\{8,18\}\\ \hline
\end{tabular}
\caption{Finiteness of $J(\{4\}^d|p)$, $2\le d\le 5$.}
\label{Ta:4dp}
\end{center}
\end{table}

\begin{table}[b]
\begin{center}
\begin{tabular}{  ||c|lc|lc|lc|lc||  } \hline
 $p $& \ &$d=2$  &\ &$d=3$   &\ &$d=4$   &\ &$d=5$  \\ \hline
 13&2,&$\emptyset$&\ &\ &\  &\ &\ &\ \\ \hline
 17&2,&\{7\}&\ &\ &\  &\ &\ &\ \\ \hline
 19&2,&\{10\}&3,&\{14\}&\ &\ &\ &\ \\ \hline
 23&2,&\{3\}&3,&$\emptyset$&4,&\{10,15\}&\ &\  \\ \hline
 29&2,&$\emptyset$&3,&$\emptyset$&4,&\{8,18\}&5,&\{18,19\} \\ \hline
 31&2,&$\emptyset$&3,&\{5,10,22,27\}&4,&\{12,18,21\}&5,&$\emptyset$ \\ \hline
 37&2,&\{18,19\} &3,&\{25,31\}&4,&\{9,18\}&5,&\{32\} \\ \hline
 41&2,&\{5,16\} &\ &\ &4,&\{17,34,35\}&5,&\{34,35,36\}\\ \hline
 43&2,&\{7,37\} &\ &\ &\ &\ &5,&$\emptyset$ \\\hline\end{tabular}
\caption{Finiteness of $J(\{5\}^d|p)$ for $2\le d\le 5$.}
\label{Ta:5dp}
\end{center}
\end{table}

\begin{table}[t]
\begin{center}
\begin{tabular}{  ||c|lc|lc|lc|lc||  } \hline
 $p $& \ &$d=2$  &\ &$d=3$   &\ &$d=4$   &\ &$d=5$  \\ \hline
 17&2,&$\emptyset$&\ &\ &\ &\ &\ &\  \\ \hline
 19&2,&\{4,7,16\}&\ &\ &\ &\ &\ &\  \\ \hline
 23&2,&$\emptyset$&3,&\{8\}&\ &\ &\ &\ \\ \hline
 29&2,&\{24\}&3,&$\emptyset$&4,&$\emptyset$&\ &\  \\ \hline
 31&2,&\{26\}&3,&\{26,27\}&4,&$\{8,13,\text{26-28}\}$&\ &\   \\ \hline
 37&2,&\{16\} &3,&$\emptyset$&4,&$\{9,15,23,33\}$&5,&$\{9,10,15,16,33,34\}$\\ \hline
 41&2,&\{28,35\}&3,&\{18\}&4,&\{32\}&5,&$\emptyset$\\ \hline
 43&2,&$\emptyset$&\ &\ &4,&$\{10,12,37,38\}$&5,&$\{10,11,14,18,29,\text{37-39}\}$\\ \hline
 47&2,&\{34\} &\ &\ &\ &\ &5,&\{36\}\\ \hline
 53&2,&\{37,49\} &\ &\ &\ &\ &5,&$\emptyset$\\ \hline\end{tabular}
\caption{Finiteness of $J(\{6\}^d|p)$, $2\le d\le 5$.}
\label{Ta:6dp}
\end{center}
\end{table}
\end{proof}

\pagebreak
\

\pagebreak

 Now we turn to non-homogeneous MZV series. We have the
following results.
\begin{prop} \label{prop:Jst}
Let $1\le s\le 10$, $2\le t\le 10$ be two positive integers. Then
there's always some prime $p>sd+2$ such that $J(s,t|p)=RJ(s,t;p)$
where $RJ(s,t)$ is the reserved set for $\zeta(s,t)$.
\end{prop}
\begin{proof}
The proposition follows from Criterion Theorem by computer
computation whose results are listed in Table~\ref{Ta:stp}. In
each entry after $(s,t)$ we list the prime $p$, the number $\gt$
used when applying Criterion Theorem to show finiteness of
$J(s,t|p)$. For these primes we found that $J(s,t|p)= RJ(s,t;p)$.

\begin{table}[h]
\begin{center}
\begin{tabular}{||c|c||c|c||c|c||c|c||c|c|| } \hline
$s,t$& $p;\gt$ &$s,t$& $p;\gt$ &$s,t$& $p;\gt$ &$s,t$& $p;\gt$
&$s,t$& $p;\gt$\\ \hline
 1,1& 13,4&2,1&13,4&3,1&11,$>$8&4,1&13,$>$8&5,1& 13,$>$8 \\ \hline
 6,1&11,$>$6 &7,1& 19,$>$6 &8,1& 13,$>$6&9,1& 17,$>$6&10,1& 31,$>$6\\ \hline
 1,2&11,4&2,2& 7,3&3,2&11,3 &4,2& 11,3&5,2& 19,4  \\ \hline
 6,2& 13,4&7,2& 19,5 &8,2& 13,5 &9,2& 23,6 &10,2& 31,6\\ \hline
 1,3& 7,5&2,3& 41,3 &3,3& 11,2 &4,3& 11,3&5,3& 11,3\\ \hline
 6,3& 13,3 &7,3&13,4&8,3& 13,4 &9,3& 13,4&10,3& 17,4\\ \hline
 1,4& 11,4 &2,4&11,3&3,4& 17,3&4,4& 11,2&5,4& 19,3\\ \hline
 6,4& 13,3 &7,4&17,3 &8,4& 17,3 &9,4& 19,3&10,4&17,3 \\ \hline
 1,5& 11,4&2,5& 11,3 &3,5&11,3 &4,5& 17,3 &5,5& 13,2\\ \hline
 6,5&47,3 &7,5& 19,3&8,5& 17,3 &9,5&17,3&10,5& 19,3\\ \hline
 1,6& 13,4&2,6& 11,3 &3,6& 17,3 &4,6& 13,3 &5,6&43,3\\ \hline
 6,6& 17,2&7,6&17,3 &8,6& 31,3 &9,6& 19,3 &10,6&19,3\\ \hline
 1,7& 17,4&2,7& 29,3 &3,7&19,3 &4,7& 37,3&5,7& 17,3\\ \hline
 6,7& 19,3 &7,7& 17,2 &8,7& 29,3&9,7& 23,3 &10,7& 23,3\\ \hline
 1,8&17,4 &2,8& 13,3 &3,8& 19,3&4,8& 23,3&5,8& 17,3\\ \hline
 6,8& 17,3&7,8& 41,3 &8,8& 19,2&9,8&29,3 &10,8& 23,3\\ \hline
 1,9& 19,4&2,9& 17,3 &3,9& 17,3 &4,9& 17,3&5,9& 17,3\\ \hline
 6,9& 19,3 &7,9& 29,3&8,9& 31,3&9,9& 29,2&10,9& 29,3\\ \hline
 1,10& 19,4 &2,10& 19,3&3,10&19,3&4,10& 31,3&5,10& 29,3\\ \hline
 6,10& 29,3&7,10& 23,3 &8,10&23,3 &9,10& 23,3 &10,10& 23,2\\ \hline
\end{tabular} \caption{Finiteness of $J(s,t|p)$ for
$2\le s,t\le 5$.} \label{Ta:stp}
\end{center}
\end{table}
\end{proof}
More data is available online \cite{online2}.
\pagebreak

\begin{prop} \label{prop:Jrst} Let $2\le r,s,t\le 5$ be
three positive integers. Then there's always some prime $p>sd+2$
such that $J(r,s,t|p)=RJ(r,s,t;p)$ where $RJ(r,s,t)$ is the
reserved set for $\zeta(r,s,t)$.
\end{prop}
\begin{proof}
The data of the proof is listed in Table~\ref{Ta:rstp}.
\begin{table}[h]
\begin{center}
\begin{tabular}{||c|c||c|c||c|c||c|c||c|c|| } \hline
 $r,s,t$& $p;\gt$ &$r,s,t$& $p;\gt$ &$r,s,t$& $p;\gt$ &$r,s,t$&
 $p;\gt$ &$r,s,t$& $p;\gt$\\ \hline
 (1,1,1)&11,$>8$ &(1,1,2)&7,6&(1,1,3)& 13,4 &(1,1,4)& 13,4 &(1,1,5)&11,4 \\ \hline
 (1,2,1)&7,$>8$ &(1,2,2)& 11,4 &(1,2,3)& 11,4 &(1,2,4)& 11,6 &(1,2,5)&11,4\\ \hline
 (1,3,1)&11,$>$6 &(1,3,2)& 11,$>$6 &(1,3,3)& 11,6 &(1,3,4)& 11,6 &(1,3,5)& 17,6\\ \hline
 (1,4,1)&13,$>$6 &(1,4,2)& 19,$>$6 &(1,4,3)& 29,$>$6 &(1,4,4)& 13,$>$6 &(1,4,5)& 19,$>$6\\ \hline
 (1,5,1)&11,$>$6 &(1,5,2)& 11,$>8$ &(1,5,3)& 13,8 &(1,5,4)& 17,$>$6 &(1,5,5)&37,$>$6\\ \hline
 (2,1,1)&7,$>$9 &(2,1,2)& 47,$>$6 &(2,1,3)& 19,$>$6 &(2,1,4)& 11,$>$5 &(2,1,5)& 13,$>$6\\ \hline
 (2,2,1)&17,6 &(2,2,2)& 11,3 &(2,2,3)& 11,3 &(2,2,4)& 13,3 &(2,2,5)& 13,3\\ \hline
 (2,3,1)&23,$>$6 &(2,3,2)& 13,4 &(2,3,3)& 13,3 &(2,3,4)& 13,3 &(2,3,5)& 13,3\\ \hline
 (2,4,1)&17,$>$6 &(2,4,2)& 17,4 &(2,4,3)& 13,4 &(2,4,4)& 29,5 &(2,4,5)& 19,4\\ \hline
 (2,5,1)&17,$>$6 &(2,5,2)& 17,4 &(2,5,3)& 13,4 &(2,5,4)& 19,5 &(2,5,5)& 17,4\\ \hline
 (3,1,1)&11,$>$6 &(3,1,2)& 11,$>$6 &(3,1,3)& 43,$>$6 &(3,1,4)& 17,$>$6 &(3,1,5)& 29,$>$6\\ \hline
 (3,2,1)&11,$>$6 &(3,2,2)& 11,3 &(3,2,3)& 13,4 &(3,2,4)& 19,3 &(3,2,5)& 17,3\\ \hline
 (3,3,1)&11,$>$6 &(3,3,2)& 13,5 &(3,3,3)& 19,4 &(3,3,4)& 23,3 &(3,3,5)& 17,3\\ \hline
 (3,4,1)&11,$>$6 &(3,4,2)& 13,4 &(3,4,3)& 13,4 &(3,4,4)& 17,3 &(3,4,5)& 19,3\\ \hline
 (3,5,1)&31,$>$6 &(3,5,2)& 17,6 &(3,5,3)& 17,4 &(3,5,4)& 17,4 &(3,5,5)& 29,4\\ \hline
 (4,1,1)&11,$>$6 &(4,1,2)& 11,$>$6 &(4,1,3)& 17,$>$6 &(4,1,4)& 13,$>$6 &(4,1,5)& 13,$>$6\\ \hline
 (4,2,1)&13,$>$6 &(4,2,2)& 23,5 &(4,2,3)& 13,4 &(4,2,4)& 37,4 &(4,2,5)& 37,5\\ \hline
 (4,3,1)&31,$>$6 &(4,3,2)& 19,5 &(4,3,3)& 17,3 &(4,3,4)& 19,4 &(4,3,5)& 17,3\\ \hline
 (4,4,1)&17,$>$6 &(4,4,2)& 23,6 &(4,4,3)& 17,4 &(4,4,4)& 19,3 &(4,4,5)& 17,3\\ \hline
 (4,5,1)&13,$>$6 &(4,5,2)& 17,5 &(4,5,3)& 23,4 &(4,5,4)& 31,3 &(4,5,5)& 17,3\\ \hline
 (5,1,1)&11,$>$6 &(5,1,2)& 13,$>$6 &(5,1,3)& 19,$>$6 &(5,1,4)& 29,$>$6 &(5,1,5)& 17,$>$6\\ \hline
 (5,2,1)&11,$>$6 &(5,2,2)& 13,5 &(5,2,3)& 17,4 &(5,2,4)& 17,4 &(5,2,5)& 17,5\\ \hline
 (5,3,1)&17,$>$6 &(5,3,2)& 17,6 &(5,3,3)& 17,4 &(5,3,4)& 17,4 &(5,3,5)& 17,4\\ \hline
 (5,4,1)&13,$>$6 &(5,4,2)& 29,6 &(5,4,3)& 17,4 &(5,4,4)& 19,3 &(5,4,5)& 17,3\\ \hline
 (5,5,1)&19,$>$6 &(5,5,2)& 17,6 &(5,5,3)& 47,4 &(5,5,4)& 23,3 &(5,5,5)& 23,3 \\
\hline
\end{tabular} \caption{Finiteness of $J(r,s,t|p)$ for
$2\le r,s,t\le 5$.} \label{Ta:rstp}
\end{center}
\end{table}
\end{proof}
More data is available online \cite{online2}. We have

Finally we make a conjecture about 2-divisible sets.
\begin{conj}
For all positive integer $d$ and $\ors\in \N^d$ the $2$-divisible
set $J(\ors|2)=\{0\}$.
\end{conj}
We put the supporting data in \cite{online3} due to its size.
We also deal with prime $3$ and $5$ separately in
\cite{online4,online5}.

\begin{prop} \label{prop:J1d}
Let $p=2,3$ or $5$. Then $J(\ors|p)$ is finite for
\begin{enumerate}
\item  $\ors=(s,t)$ where $s,t\le 20$ and $t\ge 2$.

\item  $\ors=(r,s,t)$ where $r,s,t\le 10$ and $t\ge 2$.

\item  $\ors=(q,r,s,t)$ where $q,r,s,t\le 4$ and $t\ge 2$.
\end{enumerate}
\end{prop}

\section*{Appendix I: The program code}
In this appendix we provide the pseudo code that we use to do the
computation according to our main Criterion Theorem.

Let $\ors=(s_1,\dots, s_d)$, $m=\min\{s_i: 1\le i\le d\}$, and $p$
be a fixed prime. We want to use modular arithmetic because it's
much faster than precise computation of the partial sums. So when
$n\in G_{t+1}$ we compute $p^{\wt(\ors)}H(\ors;n)$. By the
Criterion Theorem we see that if $p^{t\wt(\ors)}H(\ors;n)
\not\equiv 0 \pmod{p^{mt+m}}$ for all $n\in G_{t+1}$ then we can
take $\gt=t+1$ and we're done.

Suppose we look for $\gt$ in the range of $2$ to $e-1$ for some
positive integer $e$. In the inner for-loop of the following
pseudo code we use the formula
 $$H(s_1,\dots,s_i;n-d+i)=H(s_1,\dots,s_i;n-1-d+i)
+\frac{H(s_1,\dots,s_{i-1};n-d+(i-1))}{(n-d+i)^{s_i}}.$$

\begin{verbatim}
   Find_tau:=proc(s,p,e);

   if p<d+1 then
     t:=floor(log(d)/log(p));
   else
     t:=0;

   m:=min(s[1],...,s[d]);
   wt[i]:=s[1]+...+s[i] for i=1,...d;
   psum[i]:=0 for i=1,...d;
   ifstop=NO;

   while t<e and ifstop=NO do
     ifstop:=YES;
     psum[i]:=psum[i]*p^wt[i] for i=1,...d;
     # In the beginning of each new segment G_{t+1}, the partial sums
     # H(...; p^t-1) should be multiplied by higher powers of p.

     for n from max(d,p^t) to p^(t+1)-1 do
       for i to d do
         psum[i]:=psum[i]+p^(t*s[i])*psum[i-1]/(n-d+i)^(s[i]) mod p^(m*t+m);
       end of inner for loop;
       if psum[d]=0 then ifstop:=NO;
     end of outer for loop;
     t:=t+1;
   end while loop;

   if ifstop=YES then
     tau:=t;
   else
     e has to be larger.
\end{verbatim}

\section*{Appendix II: Conjectural structure of the set $J_1(\ors|p)$}
The following conjecture is \cite[Conj.\ 7.4]{1stpart}.
\begin{conj}\label{conj:RJ1}
Suppose $\ors\in \N^l$ such that (i) $\ors=1$, or (ii)
$\ors=(1,2,1)$, or (iii) $\ors=(2r-1,1)$  for some $r\ge 1$, or
(iv) $\ors=1^{2l}$ for some $l\ge 1$. Then $RJ(\ors)=RJ_2(\ors).$
For all other $\ors$ we have $RJ(\ors)=RJ_1(\ors).$
\end{conj}

\begin{defn} Let $\ors\in \N^l$. Define the {\em reserved density}
of the MHS $H(\ors)$ by
\begin{equation}
\text{density}(RJ(\ors); X)=\frac{\sharp\{\text{prime }p:
|\ors|+2< p<X,J(\ors|p)=RJ(\ors)\}} {\sharp\{\text{prime }p:
|\ors|+2< p<X\}}
\end{equation}
and the $m$th {\em reserved density}  by
\begin{equation}
\text{density}(RJ^m(\ors); X)=\frac{\sharp\{\text{prime }p:
|\ors|+1< p<X, \cup_{t=0}^m J_t(\ors|p)=RJ_m(\ors)\}}
{\sharp\{\text{prime }p: |\ors|+2< p<X\}}.
\end{equation}
\end{defn}
The next conjecture is \cite[Conj.\ 7.6]{1stpart}.
\begin{conj}\label{harmprimes}
Let $\ors\in \N^l$. Then
$$\operatorname{density}(RJ(\ors);\infty)=
\begin{cases}
1/\sqrt{e}, \quad &\text{ if $l=1,\ors\ge 2$},\\
1/e, &\text{ if $l=\ors=1$ or $l\ge 2$}.
\end{cases}
$$
\end{conj}

Assuming \cite[Conj.\ 7.4]{1stpart}
we have the following data to support \cite[Conj.\ 7.6]{1stpart}.
Throughout Tables \ref{Ta:den}-\ref{Ta:den3} we set $d(\ors)=\operatorname{density}(RJ(\ors);3000).$
\begin{table}[h]
\begin{center}
\begin{tabular}{||c|c|c|c|c|c|c|c|c|c|c|c||} \hline
$\ors$ & $d(\ors)$ &
$\ors$ & $d(\ors)$ &
$\ors$ & $d(\ors)$ &
$\ors$ & $d(\ors)$ &
$\ors$ & $d(\ors)$ &
$\ors$ & $d(\ors)$   \\ \hline
2 & 63.40\%

& 3& 60.37\%

& 4& 55.84\%

& 5& 61.22\%

& 6& 64.87\%

& 7& 58.78\%
\\ \hline
 8& 58.55\%

& 9& 62.06\%

& 10& 60.33\%

& 11& 63.38\%

& 12& 54.12\%

& 13& 61.65\%
\\ \hline
 14& 60.47\%

& 15& 61.41\%

& 16& 62.97\%

& 17& 60.14\%

& 18& 61.47\%

& 19& 63.83\%
\\ \hline
 20& 62.65\%

& 21& 62.18\%

& 22& 60.19\%

& 23& 63.74\%

& 24& 59.24\%

& 25& 61.14\%
\\ \hline
 26& 59.24\%

& 27& 62.32\%

& 28& 60.10\%

& 29& 63.18\%

& 30& 63.10\%

& 2,1& 36.60\%
\\ \hline
 3,1& 35.51\%

& 4,1& 37.38\%

& 5,1& 38.88\%

& 6,1& 35.36\%

& 7,1& 34.19\%

& 8,1& 34.19\%
\\ \hline
9,1& 34.74\%

& 10,1& 39.44\%

& 3,2& 39.25\%

& 4,2& 35.60\%

& 5,2& 38.17\%

&  6,2& 37.00\%
\\ \hline
7,2& 40.05\%

& 8,2& 40.38\%

& 9,2& 40.38\%

& 10,2& 39.53\%

& 4,3& 35.36\%

&  5,3& 39.11\%
\\ \hline
6,3& 36.77\%

& 7,3& 36.38\%

& 8,3& 38.73\%

& 9,3& 38.12\%

& 10,3& 39.29\%

& 5,4& 35.60\%
\\ \hline
6,4& 39.91\%

& 7,4& 42.02\%

& 8,4& 37.41\%

& 9,4& 40.94\%

& 10,4& 40.00\%

&  6,5& 37.79\%
\\ \hline
7,5& 38.35\%

& 8,5& 37.88\%

& 9,5& 34.35\%

& 10,5& 35.06\%

& 7,6& 33.41\%

&  8,6& 34.35\%
\\ \hline
9,6& 41.18\%

& 10,6& 35.38\%

& 8,7& 36.24\%

& 9,7& 36.08\%

& 10,7& 37.26\%

&  9,8& 39.62\%
\\ \hline
10,8& 42.08\%

& 10,9& 38.06\%

& 1,2& 34.96\%

& 1,3& 40.19\%

& 2,3& 36.22\%

&  1,4& 36.92\%
\\ \hline
2,4& 36.77\%

& 3,4& 38.41\%

& 1,5& 39.34\%

& 2,5& 34.66\%

& 3,5& 40.28\%

&  4,5& 40.05\%
\\ \hline
1,6& 38.41\%

& 2,6& 32.55\%

& 3,6& 39.58\%

& 4,6& 35.92\%

& 5,6& 37.56\%

&  1,7& 38.64\%
\\ \hline
2,7& 35.60\%

& 3,7& 38.50\%

& 4,7& 37.09\%

& 5,7& 35.29\%

& 6,7& 36.94\%

& 1,8& 33.02\%
\\ \hline
 2,8& 34.74\%

& 3,8& 37.09\%

& 4,8& 38.82\%

& 5,8& 39.53\%

& 6,8& 37.41\%

& 7,8& 36.94\%
\\ \hline
 1,9& 35.45\%

& 2,9& 34.98\%

& 3,9& 39.53\%

& 4,9& 36.94\%

& 5,9& 35.53\%

& 6,9& 38.12\%
\\ \hline
7,9& 34.67\%

& 8,9& 40.33\%

& 1,10& 34.27\%

& 2,10& 35.76\%

& 3,10& 38.12\%

& 4,10& 37.65\%
\\ \hline
 5,10& 35.29\%

& 6,10& 35.38\%

& 7,10& 40.09\%

& 8,10& 38.77\%

& 9,10& 33.33\%

& 1,1& 38.23\%
\\ \hline
 2,2& 39.72\%

& 3,3& 41.92\%

& 4,4& 33.26\%

& 5,5& 37.56\%

& 6,6& 37.65\%

& 7,7& 37.18\%
\\ \hline
 8,8& 36.08\%

& 9,9& 40.66\%

& 10,10& 37.35\%

& 11,11& 41.47\%

& 12,12& 35.31\%

& 13,13& 38.39\%
\\ \hline
 14,14& 44.18\%

& 15,15& 39.29\%

& 16,16& 36.90\%

& 17,17& 39.76\%

& 18,18& 35.08\%

& 19,19& 37.23\%
\\ \hline
 20,20& 38.04\%

& 21,21& 38.13\%

& 22,22& 34.53\%

& 23,23& 39.90\%

& 24,24& 35.82\%

& 25,25& 40.14\%
\\ \hline
 26,26& 40.24\%

& 27,27& 33.98\%

& 28,28& 37.11\%

& 29,29& 38.89\%

& 30,30& 35.35\%

& 1,1,1& 38.93\%
\\ \hline
1,2,1& 39.02\%

& 1,3,1& 37.38\%

& 1,4,1& 37.47\%

& 1,5,1& 37.00\%

& 1,6,1& 34.43\%

& 1,7,1& 32.32\%
\\ \hline
 1,8,1& 40.84\%

& 1,9,1& 37.09\%

& 1,10,1& 36.94\%

& 1,1,2& 36.22\%

& 1,2,2& 36.92\%

& 1,3,2& 37.47\%
\\ \hline
 1,4,2& 37.24\%

& 1,5,2& 36.77\%

& 1,6,2& 34.19\%

& 1,7,2& 34.51\%

& 1,8,2& 36.15\%

& 1,9,2& 35.53\%
\\ \hline
 1,10,2& 31.06\%

& 1,1,3& 36.45\%

& 1,2,3& 37.70\%

& 1,3,3& 35.83\%

& 1,4,3& 39.11\%

& 1,5,3& 40.28\%
\\ \hline
 1,6,3& 36.38\%

& 1,7,3& 36.62\%

& 1,8,3& 36.47\%

& 1,9,3& 39.53\%

& 1,10,3& 39.76\%

& 1,1,4& 40.75\%

\\ \hline 1,2,4& 36.77\%

& 1,3,4& 40.05\%

& 1,4,4& 38.88\%

& 1,5,4& 38.50\%

& 1,6,4& 37.32\%

& 1,7,4& 39.53\%

\\ \hline 1,8,4& 37.88\%

& 1,9,4& 38.35\%

& 1,10,4& 38.59\%

& 1,1,5& 38.88\%

& 1,2,5& 37.47\%

& 1,3,5& 34.43\%

\\ \hline 1,4,5& 40.14\%

& 1,5,5& 35.68\%

& 1,6,5& 37.65\%

& 1,7,5& 38.59\%

& 1,8,5& 34.59\%

& 1,9,5& 36.71\%

\\ \hline 1,10,5& 35.61\%

& 1,1,6& 40.98\%

& 1,2,6& 37.70\%

& 1,3,6& 36.38\%

& 1,4,6& 35.68\%

& 1,5,6& 34.12\%

\\ \hline 1,6,6& 39.76\%

& 1,7,6& 37.18\%

& 1,8,6& 36.71\%

& 1,9,6& 35.61\%

& 1,10,6& 40.09\%

& 1,1,7& 35.60\%

\\ \hline 1,2,7& 38.97\%

& 1,3,7& 39.91\%

& 1,4,7& 36.94\%

& 1,5,7& 37.88\%

& 1,6,7& 38.59\%

& 1,7,7& 36.94\%

\\ \hline 1,8,7& 37.50\%

& 1,9,7& 36.32\%

& 1,10,7& 40.19\%

& 1,1,8& 39.20\%

& 1,2,8& 39.67\%

& 1,3,8& 34.35\%

\\ \hline 1,4,8& 33.88\%

& 1,5,8& 36.71\%

& 1,6,8& 37.41\%

& 1,7,8& 38.92\%

& 1,8,8& 37.50\%

& 1,9,8& 36.88\%

\\ \hline 1,10,8& 38.77\%

& 1,1,9& 42.72\%

& 1,2,9& 36.00\%

& 1,3,9& 33.41\%

& 1,4,9& 37.18\%

& 1,5,9& 38.35\%

\\ \hline 1,6,9& 33.73\%

& 1,7,9& 39.39\%

& 1,8,9& 35.93\%

& 1,9,9& 33.57\%

& 1,10,9& 39.48\%

& 1,1,10& 35.53\%

\\ \hline 1,2,10& 36.94\%

& 1,3,10& 40.24\%

& 1,4,10& 38.12\%

& 1,5,10& 33.73\%

& 1,6,10& 37.74\%

& 1,7,10& 38.30\%

\\ \hline 1,8,10& 33.57\%

& 1,9,10& 38.06\%

& 1,10,10& 34.28\%

& 2,1,1& 46.03\%

& 2,2,1& 37.85\%

& 2,3,1& 39.34\%

\\ \hline 2,4,1& 38.64\%

& 2,5,1& 34.19\%

& 2,6,1& 35.60\%

& 2,7,1& 40.14\%

& 2,8,1& 43.19\%

& 2,9,1& 34.35\%

\\ \hline 2,10,1& 40.00\%

& 2,1,2& 34.81\%

& 2,2,2& 42.86\%

& 2,3,2& 42.16\%

& 2,4,2& 38.88\%

& 2,5,2& 40.05\%

\\ \hline 2,6,2& 37.32\%

& 2,7,2& 35.45\%

& 2,8,2& 39.06\%

& 2,9,2& 34.59\%

& 2,10,2& 33.88\%

& 2,1,3& 36.30\%

\\ \hline
 2,2,3& 38.41\%

& 2,3,3& 35.36\%

& 2,4,3& 34.43\%

& 2,5,3& 35.21\%

& 2,6,3& 40.14\%

& 2,7,3& 39.29\%
 \\ \hline
 2,8,3& 34.12\%

& 2,9,3& 33.18\%

& 2,10,3& 34.12\%

& 2,1,4& 35.36\%

& 2,2,4& 36.77\%

& 2,3,4& 40.05\%

 \\ \hline 2,4,4& 35.45\%

& 2,5,4& 37.56\%

& 2,6,4& 37.18\%

& 2,7,4& 35.76\%

& 2,8,4& 38.35\%

& 2,9,4& 34.82\%

 \\ \hline 2,10,4& 38.44\%

& 2,1,5& 37.24\%

& 2,2,5& 35.60\%

& 2,3,5& 35.92\%

& 2,4,5& 37.79\%

& 2,5,5& 40.71\%

 \\ \hline 2,6,5& 35.53\%

& 2,7,5& 34.12\%

& 2,8,5& 34.12\%

& 2,9,5& 36.32\%

& 2,10,5& 37.26\%

& 2,1,6& 36.30\%

 \\ \hline 2,2,6& 38.50\%

& 2,3,6& 37.09\%

& 2,4,6& 38.12\%

& 2,5,6& 34.35\%

& 2,6,6& 36.00\%

& 2,7,6& 37.65\%

 \\ \hline
 2,8,6& 33.73\%

& 2,9,6& 33.73\%

& 2,10,6& 36.41\%

& 2,1,7& 37.32\%

& 2,2,7& 37.79\%

& 2,3,7& 39.29\%

 \\ \hline
 \end{tabular} \caption{Density Conjecture.} \label{Ta:den}
\end{center}
\end{table}

\begin{table}[h]
\begin{center}
\begin{tabular}{||c|c|c|c|c|c|c|c|c|c|c|c||} \hline
$\ors$ & $d(\ors)$ &
$\ors$ & $d(\ors)$ &
$\ors$ & $d(\ors)$ &
$\ors$ & $d(\ors)$ &
$\ors$ & $d(\ors)$ &
$\ors$ & $d(\ors)$   \\ \hline
2,4,7& 31.06\%

& 2,5,7& 36.94\%

& 2,6,7& 36.47\%

& 2,7,7& 35.85\%

& 2,8,7& 33.02\%

& 2,9,7& 39.01\%

 \\ \hline 2,10,7& 36.64\%

& 2,1,8& 33.80\%

& 2,2,8& 36.24\%

& 2,3,8& 37.41\%

& 2,4,8& 36.47\%

& 2,5,8& 36.94\%

 \\ \hline 2,6,8& 40.80\%

& 2,7,8& 41.98\%

& 2,8,8& 32.39\%

& 2,9,8& 34.99\%

& 2,10,8& 34.99\%

& 2,1,9& 36.47\%

 \\ \hline 2,2,9& 41.65\%

& 2,3,9& 36.47\%

& 2,4,9& 38.12\%

& 2,5,9& 38.21\%

& 2,6,9& 35.38\%

& 2,7,9& 37.35\%

 \\ \hline 2,8,9& 37.12\%

& 2,9,9& 36.17\%

& 2,10,9& 34.99\%

& 2,1,10& 34.59\%

& 2,2,10& 36.94\%

& 2,3,10& 33.41\%

 \\ \hline 2,4,10& 39.39\%

& 2,5,10& 36.79\%

& 2,6,10& 38.53\%

& 2,7,10& 36.17\%

& 2,8,10& 36.88\%

& 2,9,10& 37.82\%

 \\ \hline 2,10,10& 35.07\%

& 3,1,1& 40.89\%

& 3,2,1& 40.98\%

& 3,3,1& 38.17\%

& 3,4,1& 42.39\%

& 3,5,1& 38.64\%

 \\ \hline 3,6,1& 34.27\%

& 3,7,1& 33.10\%

& 3,8,1& 35.76\%

& 3,9,1& 36.47\%

& 3,10,1& 38.12\%

& 3,1,2& 36.77\%

 \\ \hline 3,2,2& 34.90\%

& 3,3,2& 35.83\%

& 3,4,2& 35.83\%

& 3,5,2& 41.55\%

& 3,6,2& 34.51\%

& 3,7,2& 32.24\%

 \\ \hline 3,8,2& 38.35\%

& 3,9,2& 38.12\%

& 3,10,2& 35.76\%

& 3,1,3& 36.53\%

& 3,2,3& 35.83\%

& 3,3,3& 35.83\%

 \\ \hline 3,4,3& 39.20\%

& 3,5,3& 34.51\%

& 3,6,3& 37.18\%

& 3,7,3& 37.65\%

& 3,8,3& 34.59\%

& 3,9,3& 38.59\%

 \\ \hline 3,10,3& 34.67\%

& 3,1,4& 40.28\%

& 3,2,4& 41.22\%

& 3,3,4& 38.50\%

& 3,4,4& 35.45\%

& 3,5,4& 33.18\%

 \\ \hline 3,6,4& 35.76\%

& 3,7,4& 35.76\%

& 3,8,4& 36.47\%

& 3,9,4& 39.62\%

& 3,10,4& 41.27\%

& 3,1,5& 37.00\%

 \\ \hline 3,2,5& 37.09\%

& 3,3,5& 33.80\%

& 3,4,5& 36.47\%

& 3,5,5& 39.06\%

& 3,6,5& 40.47\%

& 3,7,5& 40.71\%

 \\ \hline 3,8,5& 38.92\%

& 3,9,5& 37.03\%

& 3,10,5& 37.12\%

& 3,1,6& 38.03\%

& 3,2,6& 38.26\%

& 3,3,6& 36.24\%

 \\ \hline 3,4,6& 40.71\%

& 3,5,6& 36.47\%

& 3,6,6& 35.29\%

& 3,7,6& 34.91\%

& 3,8,6& 40.80\%

& 3,9,6& 37.82\%

 \\ \hline 3,10,6& 36.41\%

& 3,1,7& 33.33\%

& 3,2,7& 37.65\%

& 3,3,7& 36.71\%

& 3,4,7& 38.82\%

& 3,5,7& 37.18\%

 \\ \hline 3,6,7& 41.98\%

& 3,7,7& 37.97\%

& 3,8,7& 34.99\%

& 3,9,7& 38.30\%

& 3,10,7& 38.06\%

& 3,1,8& 39.29\%

 \\ \hline 3,2,8& 37.18\%

& 3,3,8& 42.82\%

& 3,4,8& 34.82\%

& 3,5,8& 38.21\%

& 3,6,8& 38.21\%

& 3,7,8& 33.57\%

 \\ \hline 3,8,8& 34.99\%

& 3,9,8& 36.88\%

& 3,10,8& 34.28\%

& 3,1,9& 36.47\%

& 3,2,9& 34.59\%

& 3,3,9& 39.29\%

 \\ \hline 3,4,9& 33.73\%

& 3,5,9& 37.03\%

& 3,6,9& 36.41\%

& 3,7,9& 38.77\%

& 3,8,9& 38.30\%

& 3,9,9& 35.46\%

 \\ \hline 3,10,9& 35.07\%

& 3,1,10& 33.65\%

& 3,2,10& 35.76\%

& 3,3,10& 36.79\%

& 3,4,10& 36.08\%

& 3,5,10& 36.41\%

 \\ \hline 3,6,10& 40.19\%

& 3,7,10& 38.06\%

& 3,8,10& 35.93\%

& 3,9,10& 37.68\%

& 3,10,10& 38.39\%

& 4,1,1& 37.70\%

 \\ \hline 4,2,1& 40.52\%

& 4,3,1& 33.72\%

& 4,4,1& 38.41\%

& 4,5,1& 40.38\%

& 4,6,1& 34.51\%

& 4,7,1& 38.59\%

 \\ \hline 4,8,1& 36.47\%

& 4,9,1& 35.53\%

& 4,10,1& 36.94\%

& 4,1,2& 32.79\%

& 4,2,2& 36.77\%

& 4,3,2& 38.88\%

 \\ \hline 4,4,2& 34.98\%

& 4,5,2& 34.98\%

& 4,6,2& 34.12\%

& 4,7,2& 36.71\%

& 4,8,2& 36.24\%

& 4,9,2& 38.59\%

 \\ \hline 4,10,2& 41.04\%

& 4,1,3& 35.36\%

& 4,2,3& 38.17\%

& 4,3,3& 42.02\%

& 4,4,3& 40.14\%

& 4,5,3& 36.47\%

 \\ \hline 4,6,3& 36.00\%

& 4,7,3& 39.06\%

& 4,8,3& 37.41\%

& 4,9,3& 39.62\%

& 4,10,3& 39.86\%

& 4,1,4& 39.11\%

 \\ \hline 4,2,4& 34.74\%

& 4,3,4& 41.32\%

& 4,4,4& 35.29\%

& 4,5,4& 32.71\%

& 4,6,4& 37.88\%

& 4,7,4& 34.35\%

 \\ \hline 4,8,4& 37.26\%

& 4,9,4& 38.44\%

& 4,10,4& 37.82\%

& 4,1,5& 36.85\%

& 4,2,5& 36.15\%

& 4,3,5& $41.41^*$\%

 \\ \hline 4,4,5& 35.29\%

& 4,5,5& 34.35\%

& 4,6,5& 38.12\%

& 4,7,5& 36.56\%

& 4,8,5& 38.92\%

& 4,9,5& 38.77\%

 \\ \hline 4,10,5& 34.28\%

& 4,1,6& 34.04\%

& 4,2,6& 35.29\%

& 4,3,6& 38.35\%

& 4,4,6& 39.29\%

& 4,5,6& 37.41\%

 \\ \hline 4,6,6& 36.08\%

& 4,7,6& 34.43\%

& 4,8,6& 38.77\%

& 4,9,6& 34.52\%

& 4,10,6& 40.19\%

& 4,1,7& 36.24\%

 \\ \hline 4,2,7& 38.82\%

& 4,3,7& 37.88\%

& 4,4,7& 36.94\%

& 4,5,7& 40.57\%

& 4,6,7& 34.43\%

& 4,7,7& 38.06\%

 \\ \hline 4,8,7& 35.46\%

& 4,9,7& 37.59\%

& 4,10,7& 34.52\%

& 4,1,8& 36.71\%

& 4,2,8& 40.47\%

& 4,3,8& 37.88\%

 \\ \hline 4,4,8& 35.61\%

& 4,5,8& 39.86\%

& 4,6,8& 34.99\%

& 4,7,8& 34.28\%

& 4,8,8& 40.90\%

& 4,9,8& 32.62\%

 \\ \hline 4,10,8& 37.92\%

& 4,1,9& 36.47\%

& 4,2,9& 41.18\%

& 4,3,9& 35.85\%

& 4,4,9& 38.68\%

& 4,5,9& 34.28\%

 \\ \hline
4,6,9& 34.28\%

& 4,7,9& 36.41\%

& 4,8,9& 39.72\%

& 4,9,9& 36.26\%

& 4,10,9& 39.81\%

& 4,1,10& 37.88\%

\\ \hline  4,2,10& 40.80\%

& 4,3,10& 32.08\%

& 4,4,10& 41.61\%

& 4,5,10& 38.06\%

& 4,6,10& 37.82\%

& 4,7,10& 36.17\%

\\ \hline  4,8,10& 40.76\%

& 4,9,10& 39.57\%

& 4,10,10& 37.20\%

& 5,1,1& 36.77\%

& 5,2,1& 36.30\%

& 5,3,1& 37.94\%

\\ \hline  5,4,1& 39.67\%

& 5,5,1& 38.97\%

& 5,6,1& 38.59\%

& 5,7,1& 38.35\%

& 5,8,1& 38.35\%

& 5,9,1& 39.29\%

\\ \hline  5,10,1& 38.44\%

& 5,1,2& 37.70\%

& 5,2,2& 37.47\%

& 5,3,2& 36.38\%

& 5,4,2& 34.98\%

& 5,5,2& 38.35\%

\\ \hline  5,6,2& 36.71\%

& 5,7,2& 33.65\%

& 5,8,2& 34.59\%

& 5,9,2& 35.61\%

& 5,10,2& 36.56\%

& 5,1,3& 37.24\%

\\ \hline  5,2,3& 38.50\%

& 5,3,3& 36.15\%

& 5,4,3& 37.65\%

& 5,5,3& 34.35\%

& 5,6,3& 36.71\%

& 5,7,3& 37.65\%

\\ \hline  5,8,3& 38.44\%

& 5,9,3& 34.91\%

& 5,10,3& 38.30\%

& 5,1,4& 38.97\%

& 5,2,4& 35.92\%

& 5,3,4& $38.35^*$\%

\\ \hline  5,4,4& 41.41\%

& 5,5,4& 35.53\%

& 5,6,4& 37.18\%

& 5,7,4& 33.73\%

& 5,8,4& 37.50\%

& 5,9,4& 36.41\%

\\ \hline  5,10,4& 40.19\%

& 5,1,5& 39.20\%

& 5,2,5& 36.94\%

& 5,3,5& 39.76\%

& 5,4,5& 39.76\%

& 5,5,5& 37.88\%

\\ \hline  5,6,5& 34.91\%

& 5,7,5& 37.50\%

& 5,8,5& 39.01\%

& 5,9,5& 35.93\%

& 5,10,5& 34.75\%

& 5,1,6& 39.06\%

\\ \hline  5,2,6& 39.76\%

& 5,3,6& 35.53\%

& 5,4,6& 42.12\%

& 5,5,6& 32.78\%

& 5,6,6& 38.92\%

& 5,7,6& 39.48\%

\\ \hline  5,8,6& 37.12\%

& 5,9,6& 40.19\%

& 5,10,6& 38.30\%

& 5,1,7& 38.35\%

& 5,2,7& 35.53\%

& 5,3,7& 36.00\%

\\ \hline  5,4,7& 38.21\%

& 5,5,7& 40.09\%

& 5,6,7& 35.93\%

& 5,7,7& 32.86\%

& 5,8,7& 37.82\%

& 5,9,7& 38.53\%

\\ \hline  5,10,7& 38.39\%

& 5,1,8& 36.24\%

& 5,2,8& 37.65\%

& 5,3,8& 34.67\%

& 5,4,8& 38.21\%

& 5,5,8& 40.19\%

\\ \hline
\end{tabular} \caption{Density Conjecture.} \label{Ta:den1}
\end{center}
\end{table}

\begin{table}[h]
\begin{center}
\begin{tabular}{||c|c|c|c|c|c|c|c|c|c|c|c||} \hline
$\ors$ & $d(\ors)$ &
$\ors$ & $d(\ors)$ &
$\ors$ & $d(\ors)$ &
$\ors$ & $d(\ors)$ &
$\ors$ & $d(\ors)$ &
$\ors$ & $d(\ors)$   \\ \hline
 5,6,8& 38.06\%

& 5,7,8& 37.82\%

& 5,8,8& 35.93\%

& 5,9,8& 41.00\%

& 5,10,8& 34.60\%

& 5,1,9& 36.24\%

\\ \hline  5,2,9& 35.85\%

& 5,3,9& 34.67\%

& 5,4,9& 36.88\%

& 5,5,9& 37.12\%

& 5,6,9& 37.35\%

& 5,7,9& 39.95\%

\\ \hline  5,8,9& 39.10\%

& 5,9,9& 37.44\%

& 6,1,1& 35.60\%

& 6,2,1& 36.30\%

& 6,3,1& 39.91\%

& 6,4,1& 37.56\%

\\ \hline  6,5,1& 36.00\%

& 6,6,1& 36.94\%

& 6,7,1& 35.29\%

& 6,8,1& 39.06\%

& 6,10,1& 32.08\%

& 6,9,1& 37.50\%

\\ \hline
 6,10,3& 35.93\%

& 6,1,2& 35.83\%

& 6,2,2& 38.73\%

& 6,3,2& 38.73\%

& 6,4,2& 36.00\%

& 6,5,2& 37.41\%

\\ \hline  6,6,2& 34.59\%

& 6,7,2& 39.29\%

& 6,8,2& 36.79\%

& 6,9,2& 37.74\%

& 6,10,2& 38.53\%

& 6,1,3& 34.98\%

\\ \hline  6,2,3& 38.73\%

& 6,3,3& 35.76\%

& 6,4,3& 42.35\%

& 6,5,3& 38.12\%

& 6,6,3& 37.41\%

& 6,7,3& 32.78\%

\\ \hline  6,8,3& 39.15\%

& 6,9,3& 34.99\%

& 6,10,3& 35.93\%

& 6,1,4& 38.03\%

& 6,2,4& 34.59\%

& 6,3,4& 34.35\%

\\ \hline  6,4,4& 36.94\%

& 6,5,4& 40.71\%

& 6,6,4& 35.38\%

& 6,7,4& 37.03\%

& 6,8,4& 36.17\%

& 6,9,4& 37.35\%

\\ \hline  6,10,4& 38.77\%

& 6,1,5& 35.53\%

& 6,2,5& 35.53\%

& 6,3,5& 35.06\%

& 6,4,5& 39.29\%

& 6,5,5& 41.51\%

\\ \hline
 6,6,5& 37.50\%

& 6,7,5& 33.33\%

& 6,8,5& 37.59\%

& 6,9,5& 32.62\%

& 6,10,5& 35.46\%

& 6,1,6& 39.53\%

\\ \hline 6,2,6& 35.53\%

& 6,3,6& 40.00\%

& 6,4,6& 38.68\%

& 6,5,6& 32.55\%

& 6,6,6& 33.81\%

& 6,7,6& 40.43\%

\\ \hline 6,8,6& 34.75\%

& 6,9,6& 37.35\%

& 6,10,6& 37.44\%

& 6,1,7& 37.18\%

& 6,2,7& 38.12\%

& 6,3,7& 36.08\%

\\ \hline 6,4,7& 38.68\%

& 6,5,7& 39.01\%

& 6,6,7& 36.41\%

& 6,7,7& 30.50\%

& 6,8,7& 38.06\%

& 6,9,7& 36.73\%

\\ \hline 6,10,7& 37.92\%

& 6,1,8& 35.53\%

& 6,2,8& 38.68\%

& 6,3,8& 39.39\%

& 6,4,8& 36.41\%

& 6,5,8& 39.48\%

\\ \hline 6,6,8& 37.59\%

& 6,7,8& 38.77\%

& 6,8,8& 37.68\%

& 6,9,8& 40.76\%

& 6,10,8& 39.81\%

& 6,1,9& 38.21\%

\\ \hline 6,2,9& 37.50\%

& 6,3,9& 37.82\%

& 6,4,9& 39.24\%

& 6,5,9& 39.48\%

& 6,6,9& 35.46\%

& 6,7,9& 35.31\%

\\ \hline 6,8,9& 37.92\%

& 6,9,9& 33.41\%

& 6,10,9& 37.20\%

& 6,1,10& 38.68\%

& 6,2,10& 35.22\%

& 6,3,10& 39.01\%

\\ \hline 6,4,10& 39.95\%

& 6,5,10& 35.70\%

& 6,6,10& 36.49\%

& 6,7,10& 35.54\%

& 6,8,10& 35.07\%

& 6,9,10& 36.97\%

\\ \hline 6,10,10& 40.05\%

& 7,1,1& 42.62\%

& 7,2,1& 39.67\%

& 7,3,1& 36.38\%

& 7,4,1& 38.12\%

& 7,5,1& 39.29\%

\\ \hline 7,6,1& 38.59\%

& 7,7,1& 35.76\%

& 7,8,1& 37.97\%

& 7,9,1& 32.08\%

& 7,10,1& 36.41\%

& 7,1,2& 38.97\%

\\ \hline 7,2,2& 32.16\%

& 7,3,2& 40.94\%

& 7,4,2& 35.76\%

& 7,5,2& 33.88\%

& 7,6,2& 34.82\%

& 7,7,2& 35.14\%

\\ \hline 7,8,2& 39.86\%

& 7,9,2& 35.22\%

& 7,10,2& 34.28\%

& 7,1,3& 38.03\%

& 7,2,3& 35.29\%

& 7,3,3& 36.47\%

\\ \hline 7,4,3& 35.29\%

& 7,5,3& 34.82\%

& 7,6,3& 38.92\%

& 7,7,3& 36.32\%

& 7,8,3& 39.24\%

& 7,9,3& 35.93\%

\\ \hline 7,10,3& 40.43\%

& 7,1,4& 39.53\%

& 7,2,4& 38.82\%

& 7,3,4& 38.35\%

& 7,4,4& 37.65\%

& 7,5,4& 35.61\%

\\ \hline 7,6,4& 37.03\%

& 7,7,4& 37.59\%

& 7,8,4& 42.79\%

& 7,9,4& 36.17\%

& 7,10,4& 40.90\%

& 7,1,5& 34.59\%

\\ \hline 7,2,5& 33.41\%

& 7,3,5& 35.76\%

& 7,4,5& 36.32\%

& 7,5,5& 41.27\%

& 7,6,5& 38.30\%

& 7,7,5& 36.41\%

\\ \hline 7,8,5& 36.88\%

& 7,9,5& 40.90\%

& 7,10,5& 43.60\%

& 7,1,6& 35.76\%

& 7,2,6& 37.88\%

& 7,3,6& 43.40\%

\\ \hline 7,4,6& 36.08\%

& 7,5,6& 35.70\%

& 7,6,6& 33.81\%

& 7,7,6& 39.01\%

& 7,8,6& 36.64\%

& 7,9,6& 33.89\%

\\ \hline 7,10,6& 35.07\%

& 7,1,7& 34.12\%

& 7,2,7& 33.73\%

& 7,3,7& 39.39\%

& 7,4,7& 40.66\%

& 7,5,7& 36.64\%

\\ \hline 7,6,7& 39.48\%

& 7,7,7& 37.59\%

& 7,8,7& 41.00\%

& 7,9,7& 39.10\%

& 7,10,7& 41.23\%

& 7,1,8& 37.03\%

\\ \hline 7,2,8& 35.85\%

& 7,3,8& 35.70\%

& 7,4,8& 34.75\%

& 7,5,8& 38.06\%

& 7,6,8& 33.33\%

& 7,7,8& 37.20\%

\\ \hline 7,8,8& 37.92\%

& 7,9,8& 40.52\%

& 7,10,8& 40.28\%

& 7,1,9& 38.68\%

& 7,2,9& 33.57\%

& 7,3,9& 37.59\%

\\ \hline 7,4,9& 34.04\%

& 7,5,9& 31.44\%

& 7,6,9& 36.97\%

& 7,7,9& 33.18\%

& 7,8,9& 33.18\%

& 7,9,9& 33.18\%

\\ \hline 7,10,9& 37.68\%

& 7,1,10& 35.22\%

& 7,2,10& 37.59\%

& 7,3,10& 36.64\%

& 7,4,10& 34.75\%

& 7,5,10& 34.36\%

\\ \hline 7,6,10& 36.26\%

& 7,7,10& 35.54\%

& 7,8,10& 38.39\%

& 7,9,10& 35.54\%

& 7,10,10& 37.92\%

& 8,1,1& 35.92\%

\\ \hline 8,2,1& 38.73\%

& 8,3,1& 36.94\%

& 8,4,1& 37.65\%

& 8,5,1& 36.24\%

& 8,6,1& 34.82\%

& 8,7,1& 35.14\%

\\ \hline 8,8,1& 37.74\%

& 8,9,1& 33.81\%

& 8,10,1& 36.88\%

& 8,1,2& 36.62\%

& 8,2,2& 35.53\%

& 8,3,2& 41.18\%

\\ \hline 8,4,2& 36.24\%

& 8,5,2& 39.29\%

& 8,6,2& 37.50\%

& 8,7,2& 37.50\%

& 8,8,2& 35.93\%

& 8,9,2& 35.70\%

\\ \hline 8,10,2& 39.95\%

& 8,1,3& 41.65\%

& 8,2,3& 36.47\%

& 8,3,3& 39.76\%

& 8,4,3& 37.18\%

& 8,5,3& 37.74\%

\\ \hline 8,6,3& 37.74\%

& 8,7,3& 37.12\%

& 8,8,3& 34.99\%

& 8,9,3& 39.95\%

& 8,10,3& 37.82\%

& 8,1,4& 35.06\%

\\ \hline 8,2,4& 36.94\%

& 8,3,4& 39.29\%

& 8,4,4& 37.03\%

& 8,5,4& 37.03\%

& 8,6,4& 38.77\%

& 8,7,4& 40.19\%

\\ \hline 8,8,4& 36.88\%

& 8,9,4& 34.28\%

& 8,10,4& 42.89\%

& 8,1,5& 36.47\%

& 8,2,5& 41.18\%

& 8,3,5& 36.79\%

\\ \hline 8,4,5& 37.97\%

& 8,5,5& 36.17\%

& 8,6,5& 41.84\%

& 8,7,5& 36.17\%

& 8,8,5& 38.53\%

& 8,9,5& 38.15\%

\\ \hline 8,10,5& 36.26\%

& 8,1,6& 35.76\%

& 8,2,6& 34.43\%

& 8,3,6& 35.38\%

& 8,4,6& 32.62\%

& 8,5,6& 41.37\%

\\ \hline 8,6,6& 33.33\%

& 8,7,6& 36.17\%

& 8,8,6& 36.97\%

& 8,9,6& 38.15\%

& 8,10,6& 37.68\%

& 8,1,7& 33.73\%

\\ \hline 8,2,7& 38.21\%

& 8,3,7& 35.70\%

& 8,4,7& 33.57\%

& 8,5,7& 36.64\%

& 8,6,7& 37.59\%

& 8,7,7& 31.04\%

\\ \hline 8,8,7& 35.54\%

& 8,9,7& 34.83\%

& 8,10,7& 36.49\%

& 8,1,8& 39.86\%

& 8,2,8& 37.35\%

& 8,3,8& 39.01\%

\\ \hline 8,4,8& 33.57\%

& 8,5,8& 37.35\%

& 8,6,8& 33.41\%

& 8,7,8& 36.97\%

& 8,8,8& 39.81\%

& 8,9,8& 36.02\%

\\ \hline 8,10,8& 40.05\%

& 8,1,9& 39.72\%

& 8,2,9& 35.46\%

& 8,3,9& 35.93\%

& 8,4,9& 40.19\%

& 8,5,9& 34.83\%

\\ \hline 8,6,9& 34.60\%

& 8,7,9& 35.78\%

& 8,8,9& 36.73\%

& 8,9,9& 39.57\%

& 8,10,9& 37.92\%

& 8,1,10& 34.04\%

\\ \hline 8,2,10& 39.48\%

& 8,3,10& 34.04\%

& 8,4,10& 36.97\%

& 8,5,10& 33.65\%

& 8,6,10& 35.54\%

& 8,7,10& 36.49\%

\\ \hline
\end{tabular} \caption{Density Conjecture.} \label{Ta:den2}
\end{center}
\end{table}

\begin{table}[h]
\begin{center}
\begin{tabular}{||c|c|c|c|c|c|c|c|c|c|c|c||} \hline
$\ors$ & $d(\ors)$ &
$\ors$ & $d(\ors)$ &
$\ors$ & $d(\ors)$ &
$\ors$ & $d(\ors)$ &
$\ors$ & $d(\ors)$ &
$\ors$ & $d(\ors)$   \\ \hline
8,8,10& 36.97\%

& 8,9,10& 38.63\%

& 8,10,10& 38.48\%

& 9,1,1& 40.38\%

& 9,2,1& 36.24\%

& 9,3,1& 35.76\%

\\ \hline 9,4,1& 40.71\%

& 9,5,1& 36.71\%

& 9,6,1& 35.85\%

& 9,7,1& 33.49\%

& 9,8,1& 35.46\%

& 9,9,1& 34.28\%

\\ \hline 9,10,1& 39.95\%

& 9,1,2& 37.88\%

& 9,2,2& 37.65\%

& 9,3,2& 32.24\%

& 9,4,2& 37.41\%

& 9,5,2& 40.09\%

\\ \hline 9,6,2& 36.79\%

& 9,7,2& 35.70\%

& 9,8,2& 39.01\%

& 9,9,2& 39.72\%

& 9,10,2& 38.77\%

& 9,1,3& 40.47\%

\\ \hline 9,2,3& 41.41\%

& 9,3,3& 37.18\%

& 9,4,3& 37.50\%

& 9,5,3& 36.08\%

& 9,6,3& 37.35\%

& 9,7,3& 34.52\%

\\ \hline
9,8,3& 37.12\%

& 9,9,3& 37.59\%

& 9,10,3& 34.12\%

& 9,1,4& 34.82\%

& 9,2,4& 36.24\%

& 9,3,4& 37.50\%

\\ \hline 9,4,4& 34.43\%

& 9,5,4& 34.99\%

& 9,6,4& 37.35\%

& 9,7,4& 36.88\%

& 9,8,4& 39.01\%

& 9,9,4& 34.12\%

\\ \hline 9,10,4& 37.20\%

& 9,1,5& 35.29\%

& 9,2,5& 40.80\%

& 9,3,5& 37.26\%

& 9,4,5& 35.46\%

& 9,5,5& 39.01\%

\\ \hline
9,6,5& 40.66\%

& 9,7,5& 35.46\%

& 9,8,5& 36.49\%

& 9,9,5& 37.68\%

& 9,10,5& 39.81\%

& 9,1,6& 34.20\%

\\ \hline  9,2,6& 40.09\%

& 9,3,6& 34.75\%

& 9,4,6& 38.30\%

& 9,5,6& 35.22\%

& 9,6,6& 33.81\%

& 9,7,6& 38.15\%

\\ \hline  9,8,6& 37.20\%

& 9,9,6& 36.97\%

& 9,10,6& 39.10\%

& 9,1,7& 38.21\%

& 9,2,7& 37.35\%

& 9,3,7& 33.10\%
 \\ \hline
9,4,7& 35.70\%

& 9,5,7& 39.24\%

& 9,6,7& 35.78\%

& 9,7,7& 38.39\%

& 9,8,7& 37.20\%

& 9,9,7& 42.42\%

\\ \hline  9,10,7& 38.39\%

& 9,1,8& 35.93\%

& 9,2,8& 40.43\%

& 9,3,8& 35.46\%

& 9,4,8& 36.17\%

& 9,5,8& 36.97\%

\\ \hline  9,6,8& 38.63\%

& 9,7,8& 37.68\%

& 9,8,8& 39.57\%

& 9,9,8& 36.26\%

& 9,10,8& 33.65\%

& 9,1,9& 33.57\%

\\ \hline  9,2,9& 35.70\%

& 9,3,9& 39.24\%

& 9,4,9& 38.86\%

& 9,5,9& 40.28\%

& 9,6,9& 34.83\%

& 9,7,9& 34.60\%

\\ \hline  9,8,9& 35.31\%

& 9,9,9& 36.97\%

& 9,10,9& 39.19\%

& 9,1,10& 33.33\%

& 9,2,10& 37.82\%

& 9,3,10& 36.97\%

\\ \hline  9,4,10& 40.76\%

& 9,5,10& 38.63\%

& 9,6,10& 36.73\%

& 9,7,10& 36.02\%

& 9,8,10& 38.15\%

& 9,9,10& 37.29\%

\\ \hline  9,10,10& 36.34\%

& 10,1,1& 36.71\%

& 10,2,1& 34.12\%

& 10,3,1& 37.18\%

& 10,4,1& 32.47\%

& 10,5,1& 40.33\%

\\ \hline  10,6,1& 40.57\%

& 10,7,1& 36.64\%

& 10,8,1& 30.02\%

& 10,9,1& 39.95\%

& 10,10,1& 36.64\%

& 10,1,2& 36.47\%

\\ \hline  10,2,2& 36.71\%

& 10,3,2& 36.71\%

& 10,4,2& 42.22\%

& 10,5,2& 41.04\%

& 10,6,2& 38.30\%

& 10,7,2& 35.22\%

\\ \hline  10,8,2& 37.12\%

& 10,9,2& 38.30\%

& 10,10,2& 41.47\%

& 10,1,3& 36.47\%

& 10,2,3& 36.94\%

& 10,3,3& 39.15\%

\\ \hline  10,4,3& 37.26\%

& 10,5,3& 37.82\%

& 10,6,3& 40.90\%

& 10,7,3& 36.88\%

& 10,8,3& 35.46\%

& 10,9,3& 38.86\%

\\ \hline  10,10,3& 36.97\%

& 10,1,4& 34.82\%

& 10,2,4& 34.91\%

& 10,3,4& 39.39\%

& 10,4,4& 39.01\%

& 10,5,4& 35.22\%

\\ \hline  10,6,4& 36.17\%

& 10,7,4& 35.93\%

& 10,8,4& 37.92\%

& 10,9,4& 33.18\%

& 10,10,4& 38.39\%

& 10,1,5& 35.14\%

\\ \hline  10,2,5& 37.26\%

& 10,3,5& 34.99\%

& 10,4,5& 37.12\%

& 10,5,5& 35.70\%

& 10,6,5& 37.35\%

& 10,7,5& 35.54\%

\\ \hline  10,8,5& 34.36\%

& 10,9,5& 38.15\%

& 10,10,5& 40.76\%

& 10,1,6& 36.32\%

& 10,2,6& 42.32\%

& 10,3,6& 34.52\%

\\ \hline  10,4,6& 34.75\%

& 10,5,6& 36.17\%

& 10,6,6& 37.68\%

& 10,7,6& 39.34\%

& 10,8,6& 33.18\%

& 10,9,6& 36.73\%

\\ \hline  10,10,6& 35.78\%

& 10,1,7& 38.77\%

& 10,2,7& 37.12\%

& 10,3,7& 38.30\%

& 10,4,7& 39.24\%

& 10,5,7& 36.73\%

\\ \hline  10,6,7& 36.97\%

& 10,7,7& 34.60\%

& 10,8,7& 33.41\%

& 10,9,7& 35.54\%

& 10,10,7& 40.76\%

& 10,1,8& 31.68\%

\\ \hline  10,2,8& 35.46\%

& 10,3,8& 35.22\%

& 10,4,8& 32.23\%

& 10,5,8& 38.15\%

& 10,6,8& 36.26\%

& 10,7,8& 37.44\%

\\ \hline  10,8,8& 36.26\%

& 10,9,8& 35.78\%

& 10,10,8& 34.20\%

& 10,1,9& 39.24\%

& 10,2,9& 36.17\%

& 10,3,9& 32.46\%

\\ \hline  10,4,9& 37.44\%

& 10,5,9& 39.57\%

& 10,6,9& 37.44\%

& 10,7,9& 34.83\%

& 10,8,9& 33.89\%

& 10,9,9& 34.92\%

\\ \hline  10,10,9& 32.07\%

& 10,1,10& 37.12\%

& 10,2,10& 41.23\%

& 10,3,10& 36.26\%

& 10,4,10& 42.18\%

& 10,5,10& 35.54\%

\\ \hline  10,6,10& 39.34\%

& 10,7,10& 34.36\%

& 10,8,10& 36.34\%

& 10,9,10& 38.24\%

& 10,10,10& 40.00\%

& $\{2\}^4$& 35.60\%

\\ \hline  $\{4\}^4$& 35.85\%

& $\{6\}^4$& 33.89\%

& $\{8\}^4$& 36.67\%

& $\{10\}^4$& 37.08\%

& $\{12\}^4$& 32.45\%

& $\{14\}^4$& 41.20\%

\\ \hline  $\{16\}^4$& 34.14\%

& $\{18\}^4$& 38.05\%

& $\{20\}^4$& 36.92\%

& $\{2\}^5$& 37.79\%

& $\{4\}^5$& 39.72\%

& $\{6\}^5$& 35.71\%

\\ \hline  $\{8\}^5$& 39.95\%

& $\{10\}^5$& 41.11\%

& $\{12\}^5$& 37.29\%

& $\{14\}^5$& 33.82\%

& $\{16\}^5$& 35.70\%

& $\{18\}^5$& 37.10\%

\\ \hline  $\{20\}^5$& 35.56\%

& $\{2\}^6$& 39.76\%

& $\{4\}^6$& 39.34\%

& $\{6\}^6$& 32.70\%

& $\{8\}^6$& 40.14\%

& $\{10\}^6$& 38.02\%

\\ \hline  $\{12\}^6$& 39.27\%

& $\{14\}^6$& 40.93\%

& $\{16\}^6$& 39.90\%

& $\{18\}^6$& 36.32\%

& $\{20\}^6$& 36.16\%

& $\{2\}^7$& 36.00\%

\\ \hline  $\{4\}^7$& 39.43\%

& $\{6\}^7$& 36.21\%

& $\{8\}^7$& 35.90\%

& $\{10\}^7$& 36.98\%

& $\{2\}^8$& 42.69\%

& $\{4\}^8$& 35.00\%

\\ \hline  $\{6\}^8$& 33.89\%

& $\{8\}^8$& 39.71\%

& $\{10\}^8$& 39.36\%

& $\{12\}^8$& 35.22\%

& $\{14\}^8$& 35.91\%

& $\{16\}^8$& 36.25\%

\\ \hline  $\{18\}^8$& 32.75\%

& $\{20\}^8$& 36.80\%

& $\{2\}^9$& 42.79\%

& $\{4\}^9$& 35.08\%

& $\{6\}^9$& 33.98\%

& $\{8\}^9$& 39.76\%

\\ \hline  $\{10\}^9$& 39.31\%

& $\{12\}^9$& 35.57\%

& $\{14\}^9$& 36.00\%

& $\{16\}^9$& 36.52\%

& $\{18\}^9$& 32.57\%

& $\{20\}^9$& 36.50\%

\\ \hline  $\{2\}^{10}$& 42.79\%

& $\{4\}^{10}$& 35.17\%

& $\{6\}^{10}$& 34.14\%

& $\{8\}^{10}$& 39.61\%

& $\{10\}^{10}$& 39.51\%

& $\{12\}^{10}$& 35.66\%

\\ \hline  $\{14\}^{10}$& 36.27\%

& $\{16\}^{10}$& 36.29\%

& $\{18\}^{10}$& 32.65\%

& $\{20\}^{10}$& 36.10\%
& & & &
\\ \hline
\end{tabular} \caption{Density Conjecture.} \label{Ta:den3}
\end{center}
\end{table}

\end{document}